\documentstyle[12pt]{amsart}

\newsymbol\boxtimes 1202

\newcommand{\nc}{\newcommand}

\nc{\ch}{\operatorname{ch}}
\nc{\CHom}{{\cal{H}}om}
\nc{\Ens}{{\cal{E}}ns}
\nc{\Hom}{\operatorname{Hom}}
\nc{\opp}{\operatorname{opp}}
\nc{\parf}{\operatorname{parf}}
\nc{\rk}{\operatorname{rk}}
\nc{\Spec}{\operatorname{Spec}}
\nc{\td}{\operatorname{td}}

\nc{\BC}{{\Bbb{C}}}
\nc{\BP}{{\Bbb{P}}}
\nc{\BZ}{{\Bbb{Z}}}
\nc{\CM}{{\cal{M}}}
\nc{\CO}{{\cal{O}}}
\nc{\CT}{{\cal{T}}}
\nc{\tCM}{{\tilde{{\cal{M}}}}}
\nc{\tX}{{\tilde{X}}}

\nc{\ol}{\overline}
\nc{\ul}{\underline}
\nc{\ra}{\rightarrow}
\nc{\lra}{\longrightarrow}
\nc{\Lra}{\Longrightarrow}
\nc{\lla}{\longleftarrow}
\nc{\Llra}{\Longleftrightarrow}
\nc{\hra}{\hookrightarrow}
\nc{\iso}{\overset{\sim}{\lra}}
\nc{\rlh}{\rightleftharpoons}

\begin{document}

\title[]{A remark on virtual orientations for complete intersections}
\author{Vadim Schechtman}
\address{Max-Planck-Institut f\"{u}r Mathematik, 
Gottfried-Claren-Strasse 26, 53225 Bonn, Germany}
\email{vadik@@mpim-bonn.mpg.edu}
\maketitle


The aim of this note is to give a simple definition of genus 
zero {\em virtual orientation classes} (or {\em fundamental classes}) 
for projective complete intersections or, more generally, 
complete intersections in convex varieties, and to prove 
a {\em push forward formula} (see Lemmas 1 and 2 below) for them.

{\bf 1.} Let $X$ be a smooth complex projective variety. 
Let $X_{\beta}$ denote the stack of stable maps $f:\ C\lra X$ of genus $g$  
curves with $n$ marked points, such that $f_*([C])=\beta\in H_2(X)$. Here 
$n$ and $g$ will be fixed and suppressed from the notations. Let 
$\pi_{\beta}:\ \tX_{\beta}\lra X_{\beta}$ be the universal curve, 
and $\psi_{\beta}:\ \tX_{\beta}\lra X$ be the canonical map. 

A vector bundle $N$ over $X$ is called {\em convex} if 
$R^1\pi_{\beta*}\psi_{\beta}^*N=0$ for all $\beta$.  
In this case we will use the notation $N_{\beta}$ 
for the vector bundle $\pi_*\psi^*N$ over $X_{\beta}$. 
The variety $X$ is called convex if the tangent bundle $\CT_X$ is convex.  

{\bf Example.} Let $X=\BP^N$, and $g=0$. Then $X$ is convex, and any 
bundle of the form $N=\oplus\ \CO(l_a)\ (l_a\geq 0)$ is convex. 

{\bf 2.} Let $X$ be a convex variety and $N$ a convex bundle over $X$. 
Let $s:\ X\lra N$ be a {\em regular} (i.e. transversal to the zero section) 
section of $N$, $i:\ Y:=s^{-1}(0)\hra X$ the subscheme of its zeros.
 
Set $Y_{\beta}:=\coprod\ Y_{\gamma}$, the disjoint union over all 
$\gamma\in H_2(Y)$ mapping to $\beta$. The section $s$ induces 
the section $s_{\beta}:\ Y_{\beta}\lra N_{\beta}$, and 
we have $Y_{\beta}=s_{\beta}^{-1}(0)$. We denote by $i_{\beta}$ the 
embedding $Y_{\beta}\hra X_{\beta}$.  

{\bf Definition 1.} The {\em virtual dimension} $\dim^{virt}(Y_{\beta})$ is 
the number $\dim(X_{\beta})-\rk(N_{\beta})$. 

The section $s$ is not necessarily regular. If it is, then 
the virtual dimension of $Y_{\beta}$ coincides with its usual dimension. 

{\bf 3. Construction.} Our goal is to define certain class in the 
homology Chow group $A_{\dim^{virt}(Y_{\beta})}(Y_{\beta})$. 
We will use 
the key construction from the intersection theory, \cite{f}, 6.1. 
We have a cartesian square
$$\begin{array}{ccc}
Y_{\beta}&\overset{i_{\beta}}{\lra}&X_{\beta}\\
i_{\beta}\downarrow&\ &\downarrow s_d\\
X_{\beta}&\overset{s_{\beta}}{\lra}&N_{\beta}
\end{array}$$
where $s_0$ denotes the zero section. Let $C$ be the normal cone of $Y_{\beta}$ 
in $X_{\beta}$. We have the canonical closed embedding $C\hra i_{\beta}^*N_{\beta}$, 
whence the class $[C]\in A_d(i_{\beta}^*N)$ where $d=\dim(C)=\dim(X_{\beta})$. 
We set by definition
$$
[Y_{\beta}]^{virt}:=(p^*)^{-1}([C]).
$$
Here $p^*:\ A_{*}(Y_{\beta})\iso A_{*+\rk(N)}(i_{\beta}^*N)$ is the isomorphism 
induced by the projection $p:\ i_{\beta}^*N\lra Y_{\beta}$.

{\bf 4. Push forward formula.} {\bf Lemma 1.} {\em We have 
$$
i_{\beta*}([Y_{\beta}]^{virt})=[X_{\beta}]\cap e(N_{\beta})
$$
where $e$ denotes the Euler (top Chern) class.}

This formula was written down as a conjecture by Yu.I.Manin (talk at MPI, 
February 1997).  

{\bf Proof.} We can identify $s_0:\ X_{\beta}\lra N_{\beta}$ with the 
embedding of the scheme of zeros of the regular {\em diagonal} section 
of the vector bundle $N_{\beta}\times N_{\beta}\lra N_{\beta}$, 
and then apply \cite{f}, 6.3.4. $\Box$ 

{\bf 5. Iteration.} Let $M$ be a convex vector bundle over $Y$, 
and $j:\ Z\hra Y$ 
be the subscheme of zeros of a regular section $t:\ Y\lra M$. We set 
by definition, 
$$
\dim^{virt}(Z_{\beta}):=\dim^{virt}(Y_{\beta})-\rk(M_{\beta}). 
$$
We define
$$
[Z_{\beta}]^{virt}\in A_{\dim^{virt}(Z_{\beta})}(Z_{\beta})
$$
by the same construction as in no. 3, with $[X_{\beta}]$ 
replaced by $[Y_{\beta}]^{virt}$. 

If $M=i^*N'$ and $t$ is the restriction of a regular section 
$s':\ X\lra N'$ then $Z$ may be defined as the subscheme of zeros 
of the regular section $s\oplus s':\ X\lra N\oplus N'$, and the 
definition of $[Z_{\beta}]^{virt}$ at one step using $s\oplus s'$ coincides 
with the previous one (which used first $s$, then $t$). 

We have a generalization of Lemma 1:   

{\bf 6. Push forward formula. II.} {\bf Lemma 2.} {\em We have 
$$
j_{\beta*}([Z_{\beta}]^{virt})=[Y_{\beta}]^{virt}\cap e(M_{\beta}).
$$}

The proof is the same as in Lemma 1. 

{\bf 7.} The above constructions have obvious equivariant versions.

I am grateful to A.Goncharov and especially to Yu.I.Manin, 
for very useful discussions.

\end{document}